\documentclass[twocolumn]{autart}    

\makeatletter
\def\@biblabel#1{}
\makeatother

\usepackage{hyperref} \hypersetup{colorlinks=true, linkcolor=green, breaklinks=false, urlcolor=blue}
\usepackage{amsmath}
\usepackage{graphicx,mathrsfs,amssymb,cite,enumerate,mathrsfs,balance,float}
\usepackage{subfigure}
\usepackage{bbm}
\usepackage{dsfont}
\usepackage{color}
\usepackage[round]{natbib}
\usepackage{doi}
\usepackage[T1]{fontenc}
\usepackage{booktabs}
\usepackage{mathabx}

\newcommand{\until}[1]{\{1,\dots, #1\}}
\newcommand{\subscr}[2]{#1_{\textup{#2}}}

\newcommand{\setdef}[2]{\{#1 \; | \; #2\}}

\newcommand{\integernonnegative}{\ensuremath{\mathbb{Z}}_{\ge 0}}
\newcommand{\real}{\ensuremath{\mathbb{R}}}

\newcommand\oprocendsymbol{\hbox{$\square$}}
\newcommand\oprocend{\relax\ifmmode\else\unskip\hfill\fi\oprocendsymbol}

\renewcommand{\theenumi}{(\roman{enumi})}

\DeclareSymbolFont{bbold}{U}{bbold}{m}{n}
\DeclareSymbolFontAlphabet{\mathbbold}{bbold}

\newcommand{\Prob}{\mathbb{P}}
\newcommand{\E}{\mathbb{E}}
\newcommand{\0}{\mathbf{0}}
\newcommand{\z}{\mathbf{z}}
\newcommand{\x}{\mathbf{x}}

\usepackage[prependcaption,colorinlistoftodos]{todonotes}

\newcommand{\mcN}{\mathcal{N}}

\begin{document}

\begin{frontmatter}
\runtitle{Heterogeneous HK Model with Random Interactions}  

\title{Multidimensional Opinion Dynamics with Heterogeneous Bounded Confidences and Random Interactions\thanksref{footnoteinfo}}

\thanks[footnoteinfo]{This research is supported by the National Key Research and Development
Program of China (2022YFA1004600), the National Natural
Science Foundation of China (12288201, 72201008, 72131001, 12071465), and
the Shuanghu Laboratory (SH-2024JK31).
  Additionally, this material is based upon work supported by, or in part by, the U.~S.~Army Research Laboratory and the
  U.~S.~Army Research Office (W911NF-22-1-0233).}

\author[Cheng,Chen]{Jiangjiang Cheng}\ead{chengjiangjiang@amss.ac.cn},
\author[Chen]{Ge Chen\corauthref{cor}}\ead{chenge@amss.ac.cn}\corauth[cor]{Corresponding author.}
\author[Mei]{Wenjun Mei}\ead{mei@pku.edu.cn},
\author[Bullo]{Francesco Bullo}\ead{bullo@ucsb.edu},
\address[Cheng]{School of Mathematical Sciences, University of Chinese Academy of Sciences, Beijing 100049, China}  
\address[Chen]{Key Laboratory of Systems and Control, Academy of Mathematics and Systems Science, Chinese Academy of Sciences, Beijing 100190, China}
\address[Mei]{Department of Mechanics and Engineering Science, Peking University, Beijing, 100871, China}
\address[Bullo]{Department of Mechanical Engineering and the Center of Control, Dynamical-Systems and Computation, University of
California at Santa Barbara, CA 93106-5070, USA}


\begin{abstract}
  This paper introduces a heterogeneous multidimensional bounded confidence (BC) opinion dynamics with random pairwise interactions, whereby each pair of agents accesses each other's opinions with a specific probability.
  This revised model is motivated by the observation that the standard Hegselmann-Krause (HK) dynamics requires unrealistic all-to-all interactions at certain configurations.
  For this randomized BC opinion dynamics, regardless of initial opinions and positive confidence bounds, we show that the agents' states converge to fixed final opinions in finite time almost surely and that the convergence rate follows a negative exponential distribution in mean square.
  Furthermore, we establish sufficient conditions for the heterogeneous BC opinion dynamics with random interactions to achieve consensus in finite time.
\end{abstract}

\end{frontmatter}

\section{Introduction}

In recent decades, significant advancements have been made in agent-based opinion models, with bounded confidence (BC) models~\citep{bernardo2024bounded} standing out as one of the most compelling approaches. In this approach, the Hegselmann-Krause (HK) model~\citep{RH-UK:02} and  Deffuant-Weisbuch (DW) model~\citep{GD-DN-FA-GW:00} have gained widespread recognition.
In the HK model, each agent's opinion evolves by moving towards the arithmetic mean of the opinions within its confidence bound at each time step. The DW model assumes that a pair of agents is randomly selected at every time, and their opinions are updated if they fall within each other's confidence bounds.

The BC models have been extensively studied with a variety of methods. However, due to the interdependency and co-evolution of agents' opinions and the interagent topology over time, the theoretical analysis of the BC models remains challenging.
Most current results are restricted to the homogeneous case, whereby all agents have
the same confidence bound. The convergence and  convergence rate of
the homogeneous HK model \citep{blondel2009krause,touri2011discrete,wedin2015quadratic} and DW model have been thoroughly investigated
\citep{JL:05a,JZ-GC:15}.
Beside this body of work, many variants of the homogeneous BC model are proposed and studied for example in \cite{bernardo2022finite}, \cite{motsch2014heterophilious} and \cite{wedin2015hegselmann}.
On the other hand, convergence of the heterogeneous DW model has been proved by \cite{chen2020convergence,chen2024convergence}. For the heterogeneous HK model,
convergence properties have been established when
the confidence bound of each agent is either $0$ or $1$ \citep{BC-CW:17},  or when the amount of time each agent remains fixed is finite \citep{etesami2015game},
or when the agents' opinions are influenced by noise \citep{chen2020heterogeneous}.
However, establishing the convergence of the original heterogeneous HK model remains an open problem  \citep{mirtabatabaei2012opinion}.

It is important to acknowledge that in real-world public opinion formation
processes, it is impractical for each agent to interact with every other
agent at all times. Therefore, assuming random interactions among agents is
a more reasonable approach that has been widely pursued. For example,
\cite{fotakis2016opinion} proposes a random HK model,
where each individual randomly selected a specified number of other
individuals to interact with, and demonstrated its
convergence. \cite{granha2022opinion} explores financial market dynamics by
introducing an opinion formation model based on random interactions, which
exhibited essential real-world market characteristics both qualitatively
and quantitatively.

This paper studies the heterogeneous and multidimensional BC model with random
interactions. It should be pointed out that we do not solve the classical convergence problem of the heterogeneous HK model, but instead propose a \emph{new modification}.
In this model, each pair of agents possesses a specific
probability to access each other's opinions. We establish the finite-time
convergence of our system for any initial multidimensional opinions and
positive confidence bounds, providing the convergence rate in the
mean-square error sense. To the best of our knowledge, this study presents
the first convergence result for a heterogeneous BC model with random
interactions. Furthermore, we prove that in the final state, the opinions
of any two agents are either identical or exhibit a distance exceeding
their respective confidence bounds.  Moreover, we show that a sufficient
condition, and in certain cases a necessary condition, leads to almost sure consensus in finite time.

The paper is organized as follows. Section~\ref{Mod_sec} introduces the heterogeneous BC opinion dynamics with random interactions and presents our main convergence results. Section~\ref{pomr} contains the main proofs.
Finally, Section \ref{Conclusions} contains some concluding remarks.

\section{Our model and  main results}\label{Mod_sec}
\renewcommand{\thesection}{\arabic{section}}
In this section, we first provide preliminary knowledge and assumptions, and then propose our model and main results.

\subsection{Preliminaries and assumptions}

This paper considers the opinion evolution among a group of $n\geq 3$ agents. Let  $\mathcal{V}:=\until{n}$  denote the $n$ agents.
Assume each agent $i\in\mathcal{V}$ has a multidimensional opinion represented by a vector $\x_i(t)\in\real^d$ ($d\geq 1$) at each discrete time $t\in\integernonnegative$ \citep{parsegov2017novel}. Let $$\mathcal{S}_{\z,a}:=\setdef{\x\in\real^d}{\|\x-\z\|\leq a}$$ denote a closed sphere in $\real^d$ with center $\z\in\real^d$ and radius $a>0$, where $\|\cdot\|$ is the $\ell_2$-norm (Euclidean norm). Define $$\mathcal{S}_{\0,1}^n:=\underbrace{\mathcal{S}_{\0,1}\times \mathcal{S}_{\0,1}\times\cdots\times \mathcal{S}_{\0,1}}_n$$ to be the Cartesian product of $n$ unit spheres. Set $\x(t):=(\x_1(t),\ldots,\x_n(t))$ and assume $\x(0)\in \mathcal{S}_{\0,1}^n$. Let $r_{i}>0$ denote the \emph{confidence bound} of the agent $i$ and we assume, without loss of generality, $r_1\geq r_2\geq \cdots \geq r_n>0.$
Because $\|\x_i(0)-\x_j(0)\|\leq 2$, throughout this paper we assume that $r_n<2$, otherwise the dynamics in BC models are trivial.

The original HK model assumes that each agent can access every other agent's opinion. However, in the real world,  each agent may have no chance
 to interact with all others at every time, so the random interaction has become a rational assumption.
Also, because different agents are located in different geographical environments and have different economic and cultural backgrounds, the probability of any two agents interacting with each other may be different \citep{li2013consensus,parasnis2021convergence}. Therefore, in this paper, we assume that each pair of agents may have a different probability access to each other's opinions.
Concretely, we denote $\mathcal{G}_t=(\mathcal{V}, \mathcal{E}_t)$ as a randomly generated directed network at time $t$, where the node set $\mathcal{V}$ represents $n$ agents and the edge set $\mathcal{E}_t$ stands for the time-varying interaction relationship between agents, which is randomly selected from  $\mathcal{A}:=\setdef{(i,j)}{i,j\in\mathcal{V},i\ne {j}}$.

Next we define the probability space for our model.
Let $\mathcal{M}:=\setdef{\mathcal{E}}{\mathcal{E}\subseteq \mathcal{A}}$ denote the power set of  $\mathcal{A}$, and  $\mathcal{M}^{\infty}:=\bigtimes_{i=1}^{\infty}\mathcal{M}$ denote the infinite Cartesian product of $\mathcal{M}$.
If the initial state $\x(0)$ is deterministic, let $\Omega=\mathcal{M}^{\infty}$ be the sample space, $\mathcal{T}$ be the Tikhonov product topology on $\Omega$ \citep{knowles1973general}, $\mathcal{F}$ be the $\sigma$-algebra composed of all Borel sets on the topological space $(\Omega,\mathcal{T})$, and $\Prob$ be the probability measure on $\mathcal{F}$,
then the probability space is written as $(\Omega,\mathcal{F},\Prob)$.
If the initial state is a random variable, let $\Omega=\mathcal{S}_{\0,1}^n\times\mathcal{M}^{\infty}$ be the sample space and, similarly to the case of deterministic initial state, the probability space is defined by $(\Omega,\mathcal{F},\Prob)$.

In this paper we consider that the interaction relationship between agents is independent at different times,
and its probability has a consistent nontrivial lower bound.
In other words, we assume there exists a constant $\delta\in(0,1)$ such that
\begin{equation}\label{lowbound}
\Prob\left(\mathcal{E}_t=\mathcal{E}\right)\geq \delta, ~~\forall \mathcal{E}\subseteq \mathcal{A}, t\in\integernonnegative.
\end{equation}

We give the following example to illuminate our condition (\ref{lowbound}).
\begin{exmp}\label{examp1}
 If $\{\mathcal{G}_t=(\mathcal{V}, \mathcal{E}_t)\}$ is a sequence of i.i.d. directed Erd{\H{o}}s-R{\'e}nyi random graphs with the edge probability $p\in (0,1)$, then for any $\mathcal{E}\subseteq \mathcal{A}$ and $t\in\integernonnegative$,
  \begin{equation}\label{lowbound2}
  \begin{aligned}
    \Prob\left(\mathcal{E}_t=\mathcal{E}\right)&= \prod_{(i,j)\in\mathcal{E}} p  \prod_{(i',j')\in \mathcal{A}\setminus\mathcal{E}} (1-p)\\
    &\geq \big(\min\{p,1-p\}\big)^{n(n-1)}.
  \end{aligned}
  \end{equation}
  The condition (\ref{lowbound}) is satisfied by (\ref{lowbound2}) when we choose $\delta=\big(\min\{p,1-p\}\big)^{n(n-1)}$.
\end{exmp}

\subsection{A multidimensional and heterogeneous BC model with random interactions and main results}

Denote the neighbor set of the agent $i$ at time $t$ by
\begin{equation}\label{Neighbor}
\mcN_i(t)=\setdef{j\in\mathcal{V}}{(i,j)\in\mathcal{E}_t, \|\x_i(t)-\x_j(t)\|\leq r_i}.
\end{equation}
The dynamics of opinion evolution for the heterogeneous BC model with random interactions are as follows:
\begin{multline}\label{m1}
  \x_i(t+1)=\left(1+|\mcN_i(t)|\right)^{-1}\Big(\x_i(t)+\sum_{j\in\mcN_i(t)}\x_j(t)\Big), \\
 \forall i\in\mathcal{V}, t\geq 0,
\end{multline}
where $|\mathcal{S}|$ denotes the cardinality of set $\mathcal{S}$.
We call the system \eqref{Neighbor}-\eqref{m1} as \emph{random interaction bounded confidence (RIBC)} opinion dynamics.
If  $r_1=\cdots=r_n$ we say the RIBC model is \emph{homogeneous}; otherwise we say it is \emph{heterogenous}.

The BC models have attracted a significant amount of interest, but has shown difficult to analyze. Currently, the analysis of the BC models focuses on the homogeneous case \citep{dittmer2001consensus,wedin2015quadratic,JZ-GC:15}, while the analysis of the heterogeneous case is almost lacking \citep{mirtabatabaei2012opinion,chen2020convergence}.

The main results of this paper can be formulated as follows.

\begin{thm}{\textbf{\textup{(Convergence and convergence rate of the RIBC system)}}}
  \label{Main_result}
  Consider the RIBC system \eqref{Neighbor}-\eqref{m1} satisfying the random interaction condition \eqref{lowbound}. Then, for any initial state $\x(0)\in\mathcal{S}_{\0,1}^n$\footnotemark,
    \begin{enumerate}
    \item there exists a random variable $\x^*=(\x_1^*,\ldots,\x_n^*)\in\mathcal{S}_{\0,1}^n$ satisfying
      $\x_i^*=\x_j^*$ or $\|\x_i^*-\x_j^*\|>\max\{r_i,r_j\}$ for all $i\neq j$,
      such that $\x(t)$ converges to $\x^*$ in finite time almost surely (a.s.), and
    \item the mean-square convergence rate is negative exponential, i.e.,
   $$\sum_{i=1}^n \E\|\x_i(t)-\x_i^*\|^2 \leq 4n (1-\delta^{\lfloor T_n \rfloor})^{\lfloor \frac{t}{\lfloor T_n \rfloor+1}\rfloor},$$ where
   $\delta\in(0,1)$ is the same constant appearing in \eqref{lowbound}, and
\begin{multline*}
T_n:=\bigg[3n-2+8\left(\frac{2}{r_n}-1\right)\bigg(n+\frac{5}{3}\\
 -2\log(n+2)\bigg)\bigg]\left(\left\lceil\log_2 \frac{2}{r_n}\right\rceil+1\right)+n-2.
\end{multline*}
    \end{enumerate}
\end{thm}

\footnotetext{In fact, for any initial state $\x(0)\in\real^d$, the similar result holds. $\x(0)\in \mathcal{S}_{\0,1}^n$ is only for convenience of discussion.}

The proof of Theorem~\ref{Main_result} is postponed to Section~\ref{pomr}.
Theorem \ref{Main_result} could lead to two corollaries for consensus as follows:

\begin{cor}{\textbf{\textup{(Almost sure consensus for deterministic initial state)}}}
  \label{Cor_1}
  Consider the RIBC system \eqref{Neighbor}-\eqref{m1} satisfying the random interaction condition \eqref{lowbound}. If the largest confidence bound $r_1$ is not less than $2$, then for any initial state $\x(0)\in\mathcal{S}_{\0,1}^n$ the system reaches consensus in finite time a.s.
\end{cor}

\newcommand{\rhomin}{\subscr{\rho}{min}}

\begin{cor}{\textbf{\textup{(Almost sure consensus for random initial state)}}}
  \label{Cor_2}
  Consider the RIBC system \eqref{Neighbor}-\eqref{m1} satisfying the random interaction condition \eqref{lowbound}. Assume that the initial state $\x(0)$ is randomly distributed in $\mathcal{S}_{\0,1}^n$ and that its joint probability density has a lower bound $\rhomin>0$, that is, for any closed spheres $\mathcal{S}_{\z_i,a_i}$, $\z_i\in\real^d$, $a_i>0$ and $i\in\mathcal{V}$, with $\mathcal{S}_{\z_i,a_i}\subseteq\mathcal{S}_{\0,1}$,
  \begin{equation}\label{cor2_0}
    \Prob\Big(\bigcap_{i=1}^n \{\x_i(0)\in\mathcal{S}_{\z_i,a_i}\} \Big)
    \geq \rhomin \prod_{i=1}^n V_d(a_i),
  \end{equation}
  where $V_d(a_i)$ represents the volume of sphere with radius $a_i$ in $\real^d$\footnotemark. Then the system reaches consensus in finite time a.s. if and only if the largest confidence bound $r_1\geq 2$.
\end{cor}

\footnotetext{The volume formula of sphere in $\real^d$ with radius $a$ is $V_d(a)=\frac{\pi^{d/2}}{\Gamma(\frac{d}{2}+1)}a^d$, where $\Gamma(\frac{d}{2}+1)=\int_0^{+\infty}t^{\frac{d}{2}}\exp(-t)\mathrm{d}t$.}

\section{Proof of main results}\label{pomr}

The proof of Theorem~\ref{Main_result} requires multiple steps. We adopt the method of ``transforming the analysis of a stochastic system into the design of control algorithms'' first proposed by \cite{GC:17b}.
This method requires the construction of a new system called as the \emph{controllable interaction bounded confidence (CIBC)} opinion dynamics to help with the analysis of the RIBC system \eqref{Neighbor}-\eqref{m1}.
The CIBC system still considers the dynamics \eqref{Neighbor}-\eqref{m1}, however the underlying interaction network $\mathcal{G}_t=(\mathcal{V}, \mathcal{E}_t)$ is not randomly generated but instead treated as a control input.
In other words, we assume that the edge set $\mathcal{E}_t$ is a control input which can be chosen as arbitrary subset of $\mathcal{A}$.

The outline of the proof of Theorem~\ref{Main_result} using this method is as follows:
\begin{itemize}
  \item Firstly, Lemma~\ref{robust} provides the connection between the RIBC system and the CIBC system.
  \item Under the CIBC system,  Lemma~\ref{diameter1}  proves that two clusters of opinions with connection  can be merged into one cluster in finite time by employing suitable control inputs.
 \item Based on Lemma \ref{diameter1}, we originally propose an ingenious control algorithm (Algorithm \ref{algorithm1}) which guarantees the convergence time of the CIBC system as short as possible.
  Lemmas~\ref{timeeq}-\ref{diameter2} prove the convergence time is $O(n)$. This result and the corresponding methods are essentially different from previous works \citep{GC:17b,chen2020convergence,chen2020heterogeneous}.
  \item Finally, combining Lemma \ref{robust} with the convergence time of the CIBC system, we prove the convergence of the RIBC system, and estimate the mean square convergence rate.
\end{itemize}
It is worth pointing out that in the original heterogeneous HK model,
the opinions of agents may converge to different clusters with connection, which
makes the proof of its convergence extremely challenging.
In the CIBC system, agents with large confidence bounds can choose interaction targets to dominate the merger of clusters with connection, which leads to all agents forming different clusters without connection, and then entering  an equilibrium state quickly. With the connection between the RIBC and CIBC systems, we can prove the convergence of the RIBC system.

\subsection{Connection between the RIBC system and the CIBC system}

Given $\mathcal{S}\subseteq \mathcal{S}_{\0,1}^n$, we say $\mathcal{S}$ \emph{is reached at time $t$} if $\x(t)\in\mathcal{S}$, and \emph{is reached in the time interval $[t_1,t_2]$} if
there exists $t\in [t_1,t_2]$ such that $\x(t)\in\mathcal{S}$.

The following lemma builds a connection between the RIBC system and the CIBC system.

\begin{lem}\label{robust}{\textbf{\textup{(Connection between the RIBC system and the CIBC system)}}}
Let $\mathcal{S}\subseteq \mathcal{S}_{\0,1}^n$ be a set of states. Assume there exists a duration $t^*>0$ such that for any $\x(0)\in\mathcal{S}_{\0,1}^n$, we can find a sequence of edge sets $\mathcal{E}_0',\mathcal{E}_1',\cdots,\mathcal{E}_{t^*-1}'$ for opinion interaction which guarantees $\mathcal{S}$ is reached in the time interval $[0,t^*]$ under the CIBC system. Then, under the RIBC system \eqref{Neighbor}-\eqref{m1} with random interaction condition \eqref{lowbound}, for any initial state $\x(0)\in\mathcal{S}_{\0,1}^n$ we have
\begin{equation*}
  \Prob\left(\tau\geq t\right)\leq (1-\delta^{t^*})^{\lfloor \frac{t}{t^*+1}\rfloor},  \quad \forall t\geq 1,
\end{equation*}
where $\tau:=\min\{t':\x(t')\in\mathcal{S}\}$ is the time when $\mathcal{S}$ is first reached.
\end{lem}

Because the proof of Lemma \ref{robust} is similar to the proof of Lemma 5 in \cite{chen2020convergence}, we leave out it to save space.

According to Lemma \ref{robust}, to prove the convergence of the RIBC model \eqref{Neighbor}-\eqref{m1}, we only need to design control algorithms for the CIBC system such that a convergence set is reached.

\subsection{Design of control algorithms}
We define $\bar{n}\geq 1$ agents $\{\alpha_1,\alpha_2,\ldots,\alpha_{\bar{n}}\}$ as a \emph{cluster} at time $t$ if they have the same opinion and others have different opinions from them,
i.e., $\x_{\alpha_1}(t)=\x_{\alpha_2}(t)=\cdots=\x_{\alpha_{\bar{n}}}(t)$ and $\x_{\alpha_1}(t)\neq \x_{i}(t)$ for $i\notin \{\alpha_1,\ldots,\alpha_{\bar{n}}\}$.
Also, we say that there is a \emph{connection} between two clusters $\mathcal{C}_1,\mathcal{C}_2$  at time $t$ if the opinion distance between them is not bigger than the maximum confidence bound of the agents in them, i.e.,
$\|\x_i(t)-\x_j(t)\|\leq \max_{k\in\mathcal{C}_1\cup\mathcal{C}_2} r_k$ for $i\in \mathcal{C}_1, j\in \mathcal{C}_2$.

We first design control algorithms to merge two clusters with connection into a cluster in finite time.

\begin{lem}\label{diameter1}
 Assume $\mathcal{C}_\alpha:=\{\alpha_1,\alpha_2,\ldots,\alpha_J\}$ and $\mathcal{C}_{\beta}:=\{\beta_1,\beta_2,\ldots,\beta_K\}$ as two clusters with connection at time $t\geq 0$.
   Without loss of generality, we let $\alpha_1<\alpha_2<\cdots<\alpha_J$, $\beta_1<\beta_2<\cdots< \beta_K$, and $\alpha_1<\beta_1$.
  Set $$\mathcal{A}_{\alpha,\beta}:=\setdef{(i,j)}{i,j\in \mathcal{C}_{\alpha}\cup \mathcal{C}_{\beta}, i\ne j}.$$
  Then, under the CIBC system, there is a sequence of edge sets $\mathcal{E}_t',\mathcal{E}_{t+1}',\cdots,\mathcal{E}_{t+t^*-1}'\subseteq \mathcal{A}_{\alpha,\beta}$
  with $t^*\leq S(T+1)+1$ for opinion interaction, such that $\mathcal{C}_{\alpha}$ and $\mathcal{C}_{\beta}$ become a cluster at time $t+t^*$,
  where
  \begin{eqnarray*}
   \begin{aligned}
    &~~S:=\max\Bigg\{\Bigg\lceil \left(\frac{\|\x_{\alpha_1}(t)-\x_{\beta_1}(t)\|}{\min\{r_{\alpha_J},r_{\beta_K}\}}-1\right)\\
    &~~~~~~~~~~~~~~~~~~~~~~~~~~~~~~~~~~~~~~~~~~\times \frac{4J(K+1)}{J+K+1}+2\Bigg\rceil,0\Bigg\},\\
    &~~T:=\max\left\{\left\lceil \log_{2}\frac{\|\x_{\alpha_1}(t)-\x_{\beta_1}(t)\|}{\min\{r_{\alpha_J},r_{\beta_K}\}}\right\rceil,0\right\}.
   \end{aligned}
 \end{eqnarray*}
\end{lem}

\begin{pf}
The proof of this lemma is identical for all cases $t=0,1,2,\ldots$.  To simplify the exposition we consider only the case when $t=0$.

Firstly, we define the maximum distance between agents in the two subsets $\mathcal{C}_{\alpha}$ and $\mathcal{C}_{\beta}$ at any time $s$ as $$d(s):=\max_{i,j\in \mathcal{C}_{\alpha}\cup \mathcal{C}_{\beta}}\|\x_i(s)-\x_j(s)\|.$$
By the conditions of this lemma we have the confidence bounds $r_{\alpha_1}\geq r_{\alpha_2}\geq\cdots\geq r_{\alpha_J}$, $r_{\beta_1}\geq r_{\beta_2}\geq\cdots\geq r_{\beta_K}$, $r_{\alpha_1}\geq r_{\beta_1}$, and we have $0<d(0)\leq r_{\alpha_1}$.
To obtain our result, we only need to prove that there exists $0<t^*\leq S(T+1)+1$ such that $d(t^*)=0$.

 We let $r_{\min}^{\alpha,\beta}:=\min_{i\in \mathcal{C}_{\alpha}\cup \mathcal{C}_{\beta}} r_i.$ If $d(0)\leq r_{\min}^{\alpha,\beta}$,  we choose the set of edges $\mathcal{E}_0'=\mathcal{A}_{\alpha,\beta}$ for opinion interaction. For any $i\in \mathcal{C}_{\alpha}\cup \mathcal{C}_{\beta}$, by \eqref{Neighbor} we have $\mcN_i(0)=\mathcal{C}_{\alpha}\cup \mathcal{C}_{\beta}\setminus\{i\}$, so according to the protocol \eqref{m1}  we can get
  \begin{equation*}
    \x_i(1)=\frac{1}{J+K}\sum_{j\in \mathcal{C}_{\alpha}\cup \mathcal{C}_{\beta}}\x_j(0),
  \end{equation*}
 which indicates $d(1)=0$.

We only need to consider the case when $d(0)>r_{\min}^{\alpha,\beta}$.
Let $t_1:=\max\{\lceil \log_{2} d(0)/r_{\beta_K}\rceil,0\},$ then $t_1\leq T$.
We first show that there exists a sequence of edge sets $\mathcal{E}_{0}',\mathcal{E}_{1}',\cdots,\mathcal{E}_{t_1}'\subseteq \mathcal{A}_{\alpha,\beta}$ for opinion interaction such that
  \begin{equation}\label{x1}
d(t_1+1)<d(0)-\frac{r_{\min}^{\alpha,\beta}}{2(K+1)}.
\end{equation}
The proof of \eqref{x1} is lengthy, we put it in Appendix \ref{APP_1}.

When $J=1$, by \eqref{T+1} and \eqref{T+2} we have $d(t_1+1)=0$, and our conclusion holds. Otherwise, we let $$t_2:=\max\left\{\left\lceil \log_{2}d(t_1+1)/r_{\alpha_J}\right\rceil,0\right\}.$$ Similar to \eqref{x1} there exists a sequence of edge sets $\mathcal{E}_{t_1+1}',\mathcal{E}_{t_1+2}',\cdots,\mathcal{E}_{t_1+t_2+1}'\subseteq \mathcal{A}_{\alpha,\beta}$ with $0\leq t_2\leq T$ for opinion interaction such that
\begin{equation}\label{cg1}
  \begin{aligned}
    d(t_1+t_2+2) &<d(t_1+1)-\frac{r_{\min}^{\alpha,\beta}}{2J}\\
    &<d(0)-\frac{r_{\min}^{\alpha,\beta}}{2(K+1)}-\frac{r_{\min}^{\alpha,\beta}}{2J}.
  \end{aligned}
\end{equation}

Repeat the above processes \eqref{x1} and \eqref{cg1} until there is an integer $k^*$  satisfying
\begin{equation}\label{cg2}
    d\Big(\sum_{i=1}^{k^*}(t_i+1)\Big)\leq r_{\min}^{\alpha,\beta}.
\end{equation}
By \eqref{x1} and \eqref{cg1} we have
\begin{eqnarray*}
&&d\bigg(\sum_{i=1}^{k^*}(t_i+1)\bigg)<\\
&&~~\begin{cases}
     d(0)-\frac{k^*+1}{2}\frac{r_{\min}^{\alpha,\beta}}{2(K+1)}-\frac{k^*-1}{2}\frac{r_{\min}^{\alpha,\beta}}{2J},~~\mbox{if~$k^*$~is odd}\\
     d(0)-\frac{k^*}{2}\frac{r_{\min}^{\alpha,\beta}}{2(K+1)}-\frac{k^*}{2}\frac{r_{\min}^{\alpha,\beta}}{2J},~~~~~~~~~\mbox{if~$k^*$~is even}
\end{cases},\nonumber
\end{eqnarray*}
which indicates that
\begin{eqnarray}\label{cg3}
k^*\leq \begin{cases}
    \lceil (\frac{d(0)}{r_{\min}^{\alpha,\beta}}-1)4\frac{J(K+1)}{J+K+1}-\frac{J-K-1}{J+K+1}\rceil,~\mbox{if~$k^*$~is odd}\\
  \lceil (\frac{d(0)}{r_{\min}^{\alpha,\beta}}-1)4\frac{J(K+1)}{J+K+1}\rceil,~~~~~~~~~~~~~~~\mbox{if~$k^*$~is even}
\end{cases}
\end{eqnarray}
by \eqref{cg2}.
From \eqref{cg3} we have
\begin{equation}\label{cg4}
k^*\leq \left\lceil \left(\frac{d(0)}{r_{\min}^{\alpha,\beta}}-1\right)4\frac{J(K+1)}{J+K+1}+2\right\rceil=S.
\end{equation}

Choose $t^*=\sum_{i=1}^{k^*}(t_i+1)+1$. By \eqref{cg4} and the fact of $t_i\leq T$ for any $1\leq i\leq k^*$ we have $t^*\leq S(T+1)+1$. Also, by \eqref{cg2}
we get $d(t^*-1)\leq r_{\min}^{\alpha,\beta}$. Similar to the discussion when $d(0)\leq r_{\min}^{\alpha,\beta}$, we choose the set of edges $\mathcal{E}_{t^*-1}'=\mathcal{A}_{\alpha,\beta}$ for opinion interaction so that $d(t^*)=0$, then our conclusion holds.
\oprocend
\end{pf}

Using Lemma \ref{diameter1}, we give a control algorithm which guarantees the following event happens in finite time under the CIBC system:\\
({\bf{E1}}) The opinion values of any two agents are either same or have a distance larger than their confidence bounds.

If ({\bf{E1}}) happens, the opinions in both the RIBC system and the CIBC system will keep unchanged forever. The detailed control algorithm leading to ({\bf{E1}})
is shown in Algorithm \ref{algorithm1}.

{
\begin{table}[ht]
\refstepcounter{table}\label{algorithm1}
\begin{center}\footnotesize
\setlength{\tabcolsep}{6pt}
\begin{tabular}{p{8cm}}
\toprule
\textbf{Algorithm 1:} Control algorithm for convergence in the CIBC system \\
\midrule
\textbf{Input:}
 The number of agents $n$, the confidence bounds $r_1\geq r_2\geq\cdots\geq r_n>0$, the initial state $\x(0)=(\x_1(0),\ldots,\x_n(0))\in\mathcal{S}_{\0,1}^n$,  and the $K_0\leq n$ non-empty clusters formed by all agents at time $0$.\\
\textbf{Output:}
 The terminal time $t$.\\
\textbf{Initialization:}
 Set $t=0$, and the status of $K_0$ clusters to be \emph{non-activated}.\\
\textbf{while (E1) does not happen}\\
 \quad\textbf{Step 1:} Find two non-activated clusters $\mathcal{C}_{\alpha}$ and $\mathcal{C}_{\beta}$ satisfying: i) $\mathcal{C}_{\alpha}$ has the minimal value of $I_m(\mathcal{C}_{\alpha})$
 among all non-activated clusters which have connection with other non-activated clusters,  where $I_m(\mathcal{C}_{\alpha})=\min\setdef{i}{i\in\mathcal{C}_{\alpha}}$ denotes the minimal index of  agents in $\mathcal{C}_{\alpha}$; ii) $\mathcal{C}_{\beta}$ has the minimal value of $I_m(\mathcal{C}_{\beta})$ among all non-activated clusters  which have connection with  $\mathcal{C}_{\alpha}$.
 Then, \emph{activate} $\mathcal{C}_{\alpha}$ and $\mathcal{C}_{\beta}$, and create a \emph{merging-subset} $\mathcal{M}_{\alpha,\beta}=\mathcal{C}_{\alpha}\cup \mathcal{C}_{\beta}$.\\
 \quad\textbf{Step 2:} Repeat Step $1$ until there are no two non-activated clusters with connection.\\
 \quad\textbf{Step 3:} Using the merging process provided in the proof of Lemma \ref{diameter1},  perform a time step towards merger for all merging-subsets simultaneously.\\
 \quad\textbf{Step 4:} If a merging-subset has already finished merger, which  means a new cluster has been generated, then remove the corresponding merging-subset and set the status of the new cluster to be non-activated.\\
 \quad\textbf{Step 5:} $t\leftarrow t+1$.\\
\textbf{Until end}\\
 \bottomrule
\end{tabular}
\end{center}
\end{table}
}

For the CIBC system, the time complexity of Algorithm \ref{algorithm1} is $O(n)$\footnotemark. To prove this result, we need the following lemma.

\footnotetext{$O(n)$ means that there exists a constant $c>0$ that satisfies the terminal time $t_n\leq cn$ of Algorithm \ref{algorithm1} when $n$ is sufficiently large.}

\begin{lem}\label{timeeq}
  Consider the CIBC system. For any initial state $\x(0)=(\x_1(0),\ldots,\x_n(0))$, we execute Algorithm \ref{algorithm1} and assume that a new cluster $\mathcal{C}_{\alpha}:=\{\alpha_1,\alpha_2,\ldots,\alpha_J\}$ is generated at time $t$. Then, for the CIBC system with $J$ agents in $\mathcal{C}_{\alpha}$ whose initial state is $(\x_{\alpha_1}(0),\ldots,\x_{\alpha_J}(0))$, the terminal time for executing Algorithm \ref{algorithm1} is also $t$.
\end{lem}

\begin{pf}
  We first define the CIBC system with agents $\mathcal{V}$ and initial state $(\x_1(0),\ldots,\x_n(0))$ as System I, and the CIBC system with agents  $\mathcal{C}_{\alpha}$ and initial state $(\x_{\alpha_1}(0),\ldots,\x_{\alpha_J}(0))$
  as System II. Moreover, we define $\mathcal{C}_1,\mathcal{C}_2,\cdots,\mathcal{C}_{K_0}$ as the $K_0\leq n$ clusters formed by all agents in System I at time $0$.
By the generation process of new cluster in Algorithm \ref{algorithm1}, there exist $J_0\leq J$ clusters of the initial time which merge into the cluster  $\mathcal{C}_{\alpha}$ at time $t$ after multiple mergers.
Without loss of generality we assume $\mathcal{C}_1\cup\mathcal{C}_2\cup\cdots\cup\mathcal{C}_{J_0}=\mathcal{C}_{\alpha}$, and $\mathcal{C}_{J_0+1}\cup\mathcal{C}_{J_0+2}\cup\cdots\cup\mathcal{C}_{K_0}=\mathcal{V}\setminus\mathcal{C}_{\alpha}$.
By the definition of clusters we can obtain that at initial time the clusters formed by the agents in System II are also $\mathcal{C}_1,\mathcal{C}_2,\cdots,\mathcal{C}_{J_0}$.

Note that in Algorithm \ref{algorithm1}, the interaction between  agents within two clusters occurs only when these two clusters are merging. After the merger is finished, the opinions of all agents in the newly generated cluster remain unchanged until the next merger occurs. Therefore, all clusters can be merged but cannot be decomposed under Algorithm \ref{algorithm1}. Then, since the cluster $\mathcal{C}_{\alpha}$ is generated at time $t$ under System I and Algorithm \ref{algorithm1},  the agents in $\mathcal{C}_{\alpha}$ will not generate new clusters with the agents in  $\mathcal{V}\setminus \mathcal{C}_{\alpha}$ during the time $[0,t]$. Thus, under the System I and Algorithm \ref{algorithm1}, we assume that the agents in cluster $\mathcal{C}_{\alpha}$ generate new clusters $\mathcal{C}_1'$, $\mathcal{C}_2'$, $\ldots$, $\mathcal{C}_l'$ at time $0<t_1\leq t_2\leq \cdots \leq t_l$ with  $\mathcal{C}_l'=\mathcal{C}_{\alpha}$ and $t_l=t$.

Now we consider the System II. Due to the uniqueness of the activation order and the certainty of the merge time of clusters in Algorithm \ref{algorithm1}, we can obtain that the agents in System II will still form new clusters $\mathcal{C}_1'$, $\mathcal{C}_2'$, $\ldots$, $\mathcal{C}_l'=\mathcal{C}_{\alpha}$ at time $0<t_1\leq t_2\leq \cdots \leq t_l=t$. Because ({\bf{E1}}) must happen when all agents form a cluster, the Algorithm \ref{algorithm1} terminates at time $t$ under System II.
  \oprocend
\end{pf}

Based on Lemmas \ref{diameter1} and \ref{timeeq}, we get an upper bound of the convergence time for the CIBC system.

\begin{lem}\label{diameter2}
Consider the CIBC system. Then for any initial state, there exists a sequence of edge sets $\mathcal{E}_0',\mathcal{E}_1',\cdots,\mathcal{E}_{T-1}'$ with
\begin{multline*}
 T<T_n=\bigg[3n-2+8\left(\frac{2}{r_n}-1\right)\bigg(n+\frac{5}{3}\\
 -2\log(n+2)\bigg)\bigg]\left(\left\lceil\log_2 \frac{2}{r_n}\right\rceil+1\right)+n-2
\end{multline*}
for opinion interaction, such that \textup{({\bf{E1}})} happens.
\end{lem}

\begin{rem}
  It is worth mentioning that the best known lower and upper bound for the convergence time of the multidimensional HK model are $\Omega(n^2)$\footnotemark and $O(n^4)$, respectively \citep{bhattacharyya2013convergence,martinsson2016improved}. In our control algorithms, the upper bound for the convergence time of the CIBC system is $O(n)$, and this result is hard to be improved according to the following ``worst'' example, which also can explain why the upper bound of convergence time is $T_n$ under Algorithm \ref{algorithm1}.
\end{rem}

\footnotetext{$\Omega(n^2)$ means that there exists a constant $c>0$ that satisfies the convergence time $t_n\geq cn^2$ when $n$ is sufficiently large.}

We assume that all agents form $n$ clusters at time $0$, and a new cluster with $2$ agents is generated at time $T_1'>0$, and then a new cluster with $3$ agents is generated at time $T_2'>T_1'$, $\dots$, until the last time $T_{n-1}'>T_{n-2}'$ generates a new cluster with $n$ agents.
In this case, according to Lemma \ref{diameter1}, we have $T_1'\leq\lceil\log_2 \frac{2}{r_n}\rceil+1,$ and
\begin{multline*}
  T_i'\leq T_{i-1}'+\left\lceil\left(\frac{2}{r_n}-1\right)4\frac{2i}{i+2}+2\right\rceil\\
  \times\left(\left\lceil\log_2 \frac{2}{r_n}\right\rceil+1\right)+1, i=2,\cdots,n-1.
\end{multline*}

Therefore, from Algorithm \ref{algorithm1}, we can find a sequence of edge sets $\mathcal{E}_0',\mathcal{E}_1',\cdots,\mathcal{E}_{T-1}'$ with
\begin{multline*}
    T=T_{n-1}'\leq\bigg[1+\sum_{i=2}^{n-1}\bigg\lceil\left(\frac{2}{r_n}-1\right)\frac{8i}{i+2}+2\bigg\rceil\bigg]\\
    \times\left(\left\lceil\log_2 \frac{2}{r_n}\right\rceil+1\right)+n-2:=T_n^*
\end{multline*}
for opinion interaction, such that ({\bf{E1}}) happens.

We will prove that $T_n^*$ is an upper bound of the convergence time for all cases under Algorithm \ref{algorithm1}. On the other hand,
since the expression of $T_n^*$  seems complex, we need to make some simplification.
By the inequality $\sum_{i=1}^n \frac{1}{i}>\log(n+1)$ we have
\begin{eqnarray*}
&&\sum_{i=2}^{n-1}\bigg\lceil\left(\frac{2}{r_n}-1\right)\frac{8i}{i+2}+2\bigg\rceil\\
&&~~<\sum_{i=2}^{n-1}\left[\left(\frac{2}{r_n}-1\right)\frac{8i}{i+2}+3\right]\\\
&&~~<3(n-1)+8\left(\frac{2}{r_n}-1\right) \sum_{i=2}^{n-1}\left(1-\frac{2}{i+2}\right)\\
&&~~<3(n-1)+8\left(\frac{2}{r_n}-1\right)\left(n+\frac{5}{3}-2\log(n+2)\right),
\end{eqnarray*}
then
\begin{multline}\label{estimation}
T_n^*<\bigg[3n-2+8\left(\frac{2}{r_n}-1\right)\bigg(n+\frac{5}{3}\\
 -2\log(n+2)\bigg)\bigg]\left(\left\lceil\log_2 \frac{2}{r_n}\right\rceil+1\right)+n-2=T_n.
\end{multline}

\begin{pf*}{Proof of Lemma \ref{diameter2}}
  Consider the CIBC system. We show that for any initial state $\x(0)\in\mathcal{S}_{\0,1}^n$, if we carry out Algorithm \ref{algorithm1} then the terminal time $\tilde{T}_n$ of ({\bf{E1}}) occurrence satisfies $\tilde{T}_n\leq T_n^*$.

We adopt the induction method to prove this result. Recall that we assume $n\geq 3$. When $n=3$, the result can be obtained by Lemma \ref{diameter1} and the discussion of all situations. Next, we assume that when the number of agents $k\leq n$, the conclusion $\tilde{T}_k\leq T_k^*$ holds for any initial state, and then we prove that for $n+1$ agents, $\tilde{T}_{n+1}\leq T_{n+1}^*$ still holds for any initial state.

First, if no new cluster is generated in Algorithm \ref{algorithm1}, then $\tilde{T}_{n+1}=0<T_{n+1}^*$. We only need to consider the case when $\tilde{T}_{n+1}>0$. According to  Algorithm \ref{algorithm1}, the number of new clusters generated is finite, and the termination time is also the generation time of the last new cluster. We assume that the last new cluster is merged from $\mathcal{C}_{\alpha}$ and $\mathcal{C}_{\beta}$, with $|\mathcal{C}_{\alpha}|=J\geq 1$, $|\mathcal{C}_{\beta}|=K\geq 1$, and $J+K\leq n+1$. Let $t_{\alpha}$ and  $t_{\beta}$ are the generation times of clusters   $\mathcal{C}_{\alpha}$ and  $\mathcal{C}_{\beta}$ respectively.
By Lemma \ref{timeeq} and inductive assumption we have $t_{\alpha}=\tilde{T}_J \leq T_{J}^*$ and  $t_{\beta}=\tilde{T}_K \leq T_{K}^*$.
Then, by Algorithm \ref{algorithm1} and Lemma \ref{diameter1}, we can obtain
  \begin{equation*}
    \begin{aligned}
      &\tilde{T}_{n+1}\leq\max\{t_{\alpha},t_{\beta}\}+\bigg\lceil\left(\frac{2}{r_{n+1}}-1\right)4\frac{J(K+1)}{J+K+1}\\
      &~~~~~~~~~~~~~~~~~~~~~~~~~+2\bigg\rceil\left(\left\lceil\log_2 \frac{2}{r_{n+1}}\right\rceil+1\right)+1\\
      &~\leq\max\{T_J^*,T_K^*\}+\bigg\lceil\left(\frac{2}{r_{n+1}}-1\right)4\frac{J(K+1)}{J+K+1}\\
      &~~~~~~~~~~~~~~~~~~~~~~~~~+2\bigg\rceil\left(\left\lceil\log_2 \frac{2}{r_{n+1}}\right\rceil+1\right)+1\\
      &~\leq T_n^*+\left\lceil\left(\frac{2}{r_{n+1}}-1\right)4\frac{2n}{n+2}+2\right\rceil\\
      &~~~~~~~~~~~~~~~~~~~~~~~~~~~~~\times\left(\left\lceil\log_2 \frac{2}{r_{n+1}}\right\rceil+1\right)+1\\
      &~\leq\bigg[1+\sum_{i=2}^{n-1}\bigg\lceil\left(\frac{2}{r_{n+1}}-1\right)\frac{8i}{i+2}+2\bigg\rceil+\bigg\lceil\left(\frac{2}{r_{n+1}}-1\right)\\
      &~~~~~~~~~\frac{8n}{n+2}+2\bigg\rceil\bigg]\left(\left\lceil\log_2 \frac{2}{r_{n+1}}\right\rceil+1\right)+n-1\\
      &~=T_{n+1}^*,
    \end{aligned}
  \end{equation*}
 which indicates that our result $\tilde{T}_{n}\leq T_{n}^*$  holds by induction.

  Therefore, for any initial state, according to Algorithm \ref{algorithm1}, there exists a sequence of edge sets $\mathcal{E}_0',\mathcal{E}_1',\cdots,\mathcal{E}_{T-1}'$ with
  $T\leq T_n^*$  for opinion interaction, such that ({\bf{E1}}) happens. Combining this with \eqref{estimation} yields our result.
  \oprocend
\end{pf*}

\subsection{Final proofs}

\begin{pf*}{Proof of Theorem~\ref{Main_result}}
Let $\mathcal{S}$ be the set of all convergence points, that is,
\begin{multline*}
  \mathcal{S}:=\setdef{\x\in\mathcal{S}_{\0,1}^n}{\forall i\neq j, \x_i=\x_j\\
 \mbox{or}~ \|\x_i-\x_j\|>\max\{r_i,r_j\}},
\end{multline*}
and $\tau$ is the time when $\mathcal{S}$ is first reached under the RIBC model \eqref{Neighbor}-\eqref{m1} satisfying the random interaction condition \eqref{lowbound}.
By Lemmas ~\ref{diameter2} and \ref{robust}, we have
\begin{gather}\label{Conr_1}
\Prob\left(\tau\geq t\right)\leq (1-\delta^{\lfloor T_n\rfloor})^{\lfloor \frac{t}{\lfloor T_n\rfloor+1}\rfloor}, ~~~~\forall t\geq 1,
\end{gather}
which can be followed by $\Prob(\tau<\infty)=1$. Thus, there exists a random variable $\x^*\in\mathcal{S}$, such that $\x(t)$ converges to $\x^*$ in finite time a.s.

On the other hand, by (\ref{Conr_1}) and the total probability formula we have
\begin{align*}\label{Conr_15}
\E\|\x(t)-\x^*\|^2&=\Prob\left(\tau\geq t\right) \E\left[ \|\x(t)-\x^*\|^2 \big|\tau\geq t \right] \nonumber\\
&~~+\Prob\left(\tau<t\right) \E\left[\|\x(t)-\x^*\|^2 \big|\tau<t \right] \nonumber\\
&=\Prob\left(\tau\geq t\right) \E\left[ \|\x(t)-\x^*\|^2 \big|\tau\geq t \right] \nonumber\\
&\leq 4n (1-\delta^{\lfloor T_n\rfloor})^{\lfloor \frac{t}{\lfloor T_n\rfloor +1}\rfloor}.
\end{align*}
~~\oprocend
\end{pf*}

\begin{pf*}{Proof of Corollary~\ref{Cor_1}}
By Theorem \ref{Main_result} we have $\x(t)$ a.s. reaches a limit point $\x^*\in\mathcal{S}_{\0,1}^n$ in a finite time which satisfies either $\|\x_1^*-\x_i^*\|=0$ or $\|\x_1^*-\x_i^*\|>r_1$ for all $2\leq i\leq n$.  Because $r_1\geq 2$, we have $\|\x_1^*-\x_i^*\|=0$ for all $2\leq i\leq n$,
which indicates $\x^*$ is a consensus state.
\oprocend
\end{pf*}

\begin{pf*}{Proof of Corollary~\ref{Cor_2}}
If $r_1\geq 2$, then Corollary~\ref{Cor_1} implies that the system reaches consensus in finite time a.s.

If $r_1<2$, we choose a point $\z\in\mathcal{S}_{\0,1}$, which satisfies $\|\z\|=\frac{4+r_1}{6}$. Then equation~(\ref{cor2_0}) implies
\begin{eqnarray*}\label{cor2_1}
  \Prob\Big(\x_1(0)\in\mathcal{S}_{\z,\frac{2-r_1}{6}}, \bigcap_{i=2}^n \Big\{\x_i(0)\in\mathcal{S}_{-\z,\frac{2-r_1}{6}}\Big\} \Big)\\
  \geq \rhomin \left[V_d\left(\frac{2-r_1}{6}\right)\right]^n>0.
\end{eqnarray*}
Also, if $\x_1(0)\in\mathcal{S}_{\z,\frac{2-r_1}{6}}$ and the event $\bigcap_{i=2}^n
\{\x_i(0)\in\mathcal{S}_{-\z,\frac{2-r_1}{6}}\}$ takes place, then
$\|\x_1(0)-\x_i(0)\|\geq\frac{2+2r_1}{3}>r_1$ for $2\leq i\leq n$. In turn, this implies that the system cannot reach consensus because the agent $1$ can never interact with the agents $2,\ldots,n$.
\oprocend
\end{pf*}

\section{Conclusions}\label{Conclusions}

The BC models are a well-known class of opinion dynamics models that have attracted
significant mathematical and sociological attention in recent years.  Since
it may be unrealistic for each individual to interact with all others at
every time, this paper proposes a heterogeneous BC opinion dynamics with random
interactions, where a probability of interaction between individuals is
considered. We show that the model has finite-time convergence and give
its convergence rate in the mean-square error sense.

Numerous practical factors have not yet been considered. For instance, does
there exist an optimal probability for interaction to minimize the number
of final opinion clusters? Does there exist an optimal probability to reach
the fastest convergence rate? If we assume that each agent can only
interact with a fixed number of agents at most, what collective behavior
will the system behave? These problems are interesting and provide valuable
directions for future research.

\appendix
\section{The proof of inequality \eqref{x1}}\label{APP_1}
We consider the following two cases:

Case I: $t_1=0$. It indicates $d(0)\leq r_{\beta_K}$, then we choose the set of edges $$\mathcal{E}_{0}'=\setdef{(i,j)}{i,j\in \{\alpha_1\}\cup \mathcal{C}_{\beta}, i\ne j}$$ for opinion interaction. For any $i\in\{\alpha_1\}\cup \mathcal{C}_{\beta}$, by \eqref{Neighbor} we get $\mcN_i(0)=\{\alpha_1\}\cup \mathcal{C}_{\beta}\setminus\{i\}$, and then by \eqref{m1} we can obtain
\begin{equation}\label{T+1}
  \begin{aligned}
    \x_i(1)&=\frac{1}{K+1}\sum_{j\in\{\alpha_1\}\cup \mathcal{C}_{\beta}}\x_j(0)\\
           &=\frac{1}{K+1}(K\x_{\beta_1}(0)+\x_{\alpha_1}(0)),
  \end{aligned}
\end{equation}
while
\begin{equation}\label{notmove}
  \x_j(1)=\x_j(0),~~\forall j\in \mathcal{C}_{\alpha}\setminus \{\alpha_1\}.
\end{equation}

From \eqref{T+1} and \eqref{notmove} we can obtain that
\begin{equation*}
  \begin{aligned}
    d(1)&=\|\x_{\alpha_1}(1)-\x_{\alpha_1}(0)\|\\
        &=\left\|\frac{K}{K+1}\x_{\beta_1}(0)+\frac{1}{K+1}\x_{\alpha_1}(0)-\x_{\alpha_1}(0)\right\|\\
        &=\frac{K}{K+1}\|\x_{\beta_1}(0)-\x_{\alpha_1}(0)\|\\
        &=d(0)-\frac{d(0)}{K+1}<d(0)-\frac{r_{\min}^{\alpha,\beta}}{2(K+1)},
  \end{aligned}
\end{equation*}
which indicates \eqref{x1}.

Case II: $t_1\geq 1$. By the definition of $t_1$ we have
\begin{multline}\label{t1}
 t_1-1<\log_2 \frac{d(0)}{r_{\beta_K}} \leq t_1 \\
  \iff 2^{t_1-1}r_{\beta_K}<d(0)\leq 2^{t_1}r_{\beta_K}.
\end{multline}
We design the sequence of edge sets $$\mathcal{E}_0'=\mathcal{E}_1'=\cdots=\mathcal{E}_{t_1-1}'=\left\{(\alpha_1,\beta_K),(\beta_K,\alpha_1)\right\}$$ and $$\mathcal{E}_{t_1}'=\setdef{(i,j)}{i,j\in \{\alpha_1\}\cup \mathcal{C}_{\beta}, i\ne j}$$ for opinion interaction at times $0,1,\ldots,t_1$.

Using (\ref{t1}) and the protocol \eqref{Neighbor}-\eqref{m1} repeatedly we obtain
\begin{eqnarray}\label{s1}
\left\{
\begin{array}{ll}
\x_i(s)=\x_i(0), ~~\forall i\in \mathcal{C}_{\beta} \\
\x_{\alpha_1}(s)=(1-1/2^s)\x_{\beta_1}(0)+\x_{\alpha_1}(0)/2^s
\end{array}
\right.
\end{eqnarray}
for $s=1,\ldots,t_1$, and
\begin{eqnarray}\label{T+2}
  \x_i(t_1+1)&&=\x_{\alpha_1}(t_1+1)\\
  &&=\frac{1}{K+1}(K\x_{\beta_1}(0)+\x_{\alpha_1}(t_1)), \forall i\in \mathcal{C}_{\beta}.\nonumber
\end{eqnarray}
Besides, the agents' opinions in $\mathcal{C}_{\alpha}\setminus\{\alpha_1\}$ remain unchanged at time $0,\ldots,t_1+1$.

By \eqref{t1} and \eqref{s1} we have
 $$\|\x_{\alpha_1}(t_1)-\x_{\beta_1}(0)\|=\frac{1}{2}\|\x_{\alpha_1}(t_1-1)-\x_{\beta_1}(0)\|>\frac{r_{\beta_K}}{2},$$ and then by \eqref{T+2} we have
\begin{equation*}
  \begin{aligned}
    d(t_1&+1)=\|\x_{\alpha_1}(t_1+1)-\x_{\alpha_1}(0)\|\\
         &=\Big\|\frac{K}{K+1}\x_{\beta_1}(0)+\frac{1}{K+1}\x_{\alpha_1}(t_1)-\x_{\alpha_1}(0)\Big\|\\
         &=\Big\|\x_{\beta_1}(0)-\x_{\alpha_1}(0)-\frac{1}{K+1}(\x_{\beta_1}(0)-\x_{\alpha_1}(t_1))\Big\|\\
         &=\|\x_{\beta_1}(0)-\x_{\alpha_1}(0)\|-\frac{1}{K+1}\|\x_{\beta_1}(0)-\x_{\alpha_1}(t_1)\|\\
         &<d(0)-\frac{1}{2(K+1)}r_{\beta_K}\leq d(0)-\frac{r_{\min}^{\alpha,\beta}}{2(K+1)},
  \end{aligned}
\end{equation*}
where the fourth equality use the fact that the two vectors with difference are collinear. This equation also indicates \eqref{x1}.
\qquad\hfill

\bibliographystyle{plainnat}
\bibliography{alias,FB,Main}

\end{document}